# Symmetrizers and Continuity of Stable Subspaces for Parabolic-Hyperbolic Boundary Value Problems

Guy Métivier,[*] Kevin Zumbrun[†]

**Abstract**

In this paper we prove the continuity of stable subspaces associated to parabolic-hyperbolic boundary value problems, for limiting values of parameters. The analysis is based on the construction performed in [MZ] of Kreiss' type symmetrizers.

## 1 Introduction

This paper is motivated by the analysis of multidimensional small viscosity boundary layers or shock profiles (see [MZ], [GMWZ1] [GMWZ2]). The stability criterium for layers or profiles is described through the nonvanishing of an Evans function. After freezing the tangential coefficients, performing a tangential Fourier-Laplace transform and using a reduction to a normal form, the problem is reduced to an ordinary $N \times N$ constant coefficient differential system in the normal variable:

$$(1.1) \qquad \partial_x U = \mathcal{G}(p, \zeta)U + F \quad \text{for} \quad x > 0,$$

where $p$ are parameters and $\zeta$ denote the tangential frequency variables (see [MZ]). Typically, $\zeta = (\tau - i\gamma, \eta) \in \mathbb{C} \times \mathbb{R}^{d-1}$ where $\tau - i\gamma$ is the Fourier-Laplace frequency dual to time $t$ and $\eta = (\eta_1, \ldots, \eta_{d-1})$ are the Fourier frequencies dual to the spatial tangential variables $y$. The equation is supplemented by boundary conditions

$$(1.2) \qquad \Gamma(p, \zeta)U_{|x=0} = 0.$$

[*]MAB Université de Bordeaux I, 33405 Talence Cedex, France; metivier@math.u-bordeaux.fr

[†]Indiana University, Bloomington, IN 47405; kzumbrun@indiana.edu. Research partially supported under NSF grant number DMS-0070765.

For $\rho = |\zeta| > 0$ and $\gamma \geq 0$, $\mathcal{G}(p,\zeta)$ has no purely imaginary eigenvalues and the bounded solutions of the homogeneous equation (1.1) are

$$U(x) = e^{x\mathcal{G}(p,\zeta)}U(0), \qquad U(0) \in \mathbb{E}_-(p,\zeta)$$

where $\mathbb{E}_-(p,\zeta)$ is the stable subspace at $+\infty$ of $\mathcal{G}(p,\zeta)$, that is the space generated by the generalized eigenvectors associated to eigenvalues belonging the the half space $\{\operatorname{Re} \mu < 0\}$. The stability condition, that is the well posedness of the boundary value problem (1.1) (1.2) in $L^p$ spaces, reads

(1.3) $$\mathbb{C}^N = \mathbb{E}_-(p,\zeta) \oplus \ker \Gamma(p,\zeta),$$

or

(1.4) $$D(p,\zeta) := \det\bigl(\mathbb{E}_-(p,\zeta), \ker \Gamma(p,\zeta)\bigr) \neq 0$$

where, for linear subspaces $E$ and $F$ of $\mathbb{C}^N$, $\det(E,F)$ is zero if $\dim E + \dim F \neq N$ and equal to the determinant obtained by taking orthonormal basis in $E$ and $F$ when $\dim E + \dim F = N$; this is *necessary* for multidimensional stability. The function $D$ is commonly known as an *Evans function*.

The *uniform stability condition* requires that $D$ (or a rescaled version for high frequencies) be uniformly bounded from below by a positive constant when $\rho > 0$ and $\gamma \geq 0$. This has been shown in several contexts to be *sufficient* for multidimensional stability; see [Z], [MZ],[GMWZ1–2].

Of particular importance is the analysis for *low* frequencies where the hyperbolic character becomes crucial. There one uses polar coordinates

(1.5) $$\zeta = \rho\check{\zeta}, \quad \rho = |\zeta|.$$

The main result of this paper is:

**Theorem 1.1.** *Under Assumption 3.1 below, the linear bundle $\check{\mathbb{E}}(p,\check{\zeta},\rho) := \mathbb{E}_-(p,\rho\check{\zeta})$ has a continuous extension to $\rho = 0$, $\check{\gamma} \geq 0$.*

This result was known under additional assumptions such as strict hyperbolicity of the hyperbolic part or constant rank of branching points, see [Rou], [ZS], or [Z]. The Assumptions 3.1 are the somewhat minimal conditions of hyperbolicity and parabolicity as in [MZ]. An important corollary is that

(1.6) $$D(p,\rho\check{\zeta}) = \beta\Delta(p,\check{\zeta}) + o(1), \qquad \text{as} \quad \rho \to 0$$

uniformly for $|\check\zeta| = 1$ with $\check\gamma \geq 0$, where $\Delta(p, \check\zeta)$ is the Lopatinski determinant of the appropriate limiting hyperbolic boundary value problem, and $\beta$ is a transversality coefficient for the inner layer ordinary differential equation (the viscous shock or boundary profile problem); indeed, the two statements are essentially equivalent. This reduces the somewhat abstract uniform stability condition, at least in the crucial low frequency regime, to a pair of checkable conditions, associated respectively with outer (hyperbolic) and inner (stationary parabolic) problems in the formal matched asymptotic approximation of the boundary layer, and is one of the central conclusions of [ZS],[Z],[Rou].

In the original planar shock stability analysis of [Z], which was carried out by explicit estimates on the Laplace–Fourier inversion formula, continuous extension of stable subspaces was also used in an essential way in obtaining the linearized decay estimates from which nonlinear stability ultimately follows. In the more recent curved boundary layer analysis of [MZ], carried out by Kreiss symmetrizer techniques, it was shown that the uniform stability condition alone suffices, whether or not continuity holds, yielding applications to much more general situations than considered in [Z]; this was possible because of the flexibility of the symmetrizer method. This abstract result, however, sidestepped the issue of checkability of the uniform stability condition/validity of the fundamental relation (1.6).

With the establishment of Theorem 1.1, we now resolve this issue completely, under the same minimal hypotheses considered in [MZ]. The matrix perturbation problem associated with determination of stable subspaces is a rather complicated one involving many parameters, and would be difficult if not impossible to describe completely in this generality. In the inviscid case, this problem can be sidestepped thanks to the block structure condition enjoyed by systems with constant multiplicity characteristics; see [Mé3]. However, the additional parameter $\rho$ generically destroys this block structure, and so the issue apparently arises again. A notable aspect of [MZ] was the construction of Kreiss symmetrizers *without the use of block structure* in either $\gamma$ or $\rho$: that is, the identification of different roles played by Laplace and Fourier frequencies in the analysis. Here, we show that, somewhat surprisingly, the same symmetrizer construction used in [MZ] to avoid the need for continuity of subspaces in fact establishes that continuity! That is, the structure of stable subspaces is already hidden in the details of the nonstandard symmetrizer construction.

As the proof of the main theorem relies on the use of symmetrizers for $\mathcal{G}(p, \zeta)$, we first review the construction performed in [MZ] of Kreiss' type symmetrizers for parabolic-hyperbolic boundary value problems and we prove the new continuity theorem next.

## 2 Symmetrizers

Recall first the essence of the "method of symmetrizers" as it applies to general boundary value problems on the half line $\{x \geq 0\}$:

(2.1) $$\partial_x u = G(x)u + f, \quad \Gamma u(0) = 0.$$

Here, $u$ and $f$ are functions on $[0, \infty[$ values in some Hilbert space $\mathcal{H}$, and $G(x)$ is a $C^1$ family of (possibly unbounded) operators defined on $\mathcal{D}$, dense subspace of $\mathcal{H}$.

A *symmetrizer* is a family of $C^1$ functions $x \mapsto S(x)$ with values in the space of operators in $\mathcal{H}$ such that there are $C_0$, $\lambda > 0$, $\delta > 0$ and $C_1$ such that

(2.2) $$\forall x, \quad S(x) = S(x)^* \quad \text{and} \quad |S(x)| \leq C_0,$$

(2.3) $$\forall x, \quad 2\operatorname{Re} S(x)G(x) + \partial_x S(x) \geq 2\lambda \operatorname{Id},$$

(2.4) $$S(0) \geq \delta \operatorname{Id} - C_1 \Gamma^* \Gamma.$$

In (2.2), the norm of $S(x)$ is the norm in the space of bounded operators in $\mathcal{H}$. Similarly $S^*(x)$ is the adjoint operator of $S(x)$. The notation $\operatorname{Re} T = \frac{1}{2}(T + T^*)$ is used in (2.3) for the real part of an operator $T$. When $T$ is unbounded, the meaning of $\operatorname{Re} T \geq \lambda$, is that all $u \in \mathcal{D}$ belongs to the domain of $T$ and satisfies

(2.5) $$\operatorname{Re}(Tu, u) \geq \lambda |u|^2.$$

The property (2.3) has to be understood in this sense.

**Lemma 2.1.** *If there is a symmetrizer $S$, then for all $u \in C_0^1([0, \infty[; \mathcal{H}) \cap C^0([0, \infty[; \mathcal{D})$, one has*

(2.6) $$\lambda \|u\|^2 + \delta |u(0)|^2 \leq \frac{C_0^2}{\lambda} \|\partial_x u - Gu\|^2 + C_1 |\Gamma u(0)|^2.$$

Here, $|\cdot|$ is the norm in $\mathcal{H}$ and $\|\cdot\|$ the norm in $L^2([0, \infty[; \mathcal{H})$.

*Proof.* Taking the scalar product of $Su$ with the equation (2.1) and integrating over $[0, \infty[$, (2.2) implies

(2.7) $$-(S(0)u(0), u(0)) = \int \partial_x (Su, u) dx$$
$$= \int ((2\operatorname{Re} SG + \partial_x S)u, u) dx + 2\operatorname{Re} \int (Sf, u) dx.$$

By (2.3),
$$\int ((2\operatorname{Re} SG + \partial_x S)u, u)dx \geq 2\lambda \|u\|^2.$$

By (2.4) and the boundary condition $\Gamma u(0) = 0$,
$$(S(0)u(0), u(0)) \geq \delta |u(0)|^2 - C_1|\Gamma u(0)|^2.$$

By (2.2)
$$2\left|\int (Sf, u)dx\right| \leq 2C_0\|f\|\,\|u\| \leq \frac{C_0^2}{\lambda}\|f\|^2 + \lambda\|u\|^2.$$

Thus the identity (2.7) implies the energy estimate (2.6). $\square$

In applications to evolution boundary value problems, $G(x)$ is a tangential differential system $G(t, y, x, \partial_t, \partial_y)$. One usually studies (2.1) in weighted spaces $e^{\gamma t}L^2$, with $\gamma > 0$ for forward propagation. With $u_\gamma = e^{-\gamma t}u$ and $f_\gamma = e^{-\gamma t}f$, the equation (2.1) is transformed into

(2.8) $$\partial_x u_\gamma = G_\gamma(x)u_\gamma + f_\gamma, \quad \Gamma u(0) = 0.$$

with $G_\gamma(t, y, x, \partial_t, \partial_y) = G(t, y, x, \partial_t + \gamma, \partial_y)$, which is studied in $L^2$. When $G$ has constant coefficients in $(t, y)$, the equation (2.8) transformed into a system of o.d.e. in $x$ using a tangential space time Fourier transform. For "slowly" varying coefficients, one first freeze the coefficients and then perform the Fourier transform, yielding a family of systems depending on the parameters $(t, z)$ and the frequencies $\zeta = (\tau, \eta, \gamma)$:

(2.9) $$\partial_x \widehat{u}_\gamma = \widehat{G}(t, y, x, \zeta)\widehat{u}_\gamma + \widehat{f}_\gamma = 0, \quad \Gamma \widehat{u}_\gamma = 0,$$

where $\widehat{G}(t, y, x, \zeta) = G(t, y, x, i\tau + \gamma, i\eta)$.

In this framework, one can look for symmetrizers for the system (2.9). This is a matrix $\widehat{S}(t, y, x, \zeta)$ such that

(2.10)
$$\widehat{S}(t, y, x, \zeta) = \widehat{S}(t, y, x, \zeta) \quad \text{and} \quad |S(t, y, x, \zeta)| \leq C_0,$$
$$\forall \gamma \geq 0, \quad 2\operatorname{Re}\widehat{S}(t, y, x, \zeta)\widehat{G}(t, y, x\zeta) + \partial_x \widehat{S}(t, y, x, \zeta) \geq \gamma \operatorname{Id},$$
$$\widehat{S}(t, y, 0, \zeta) \geq \delta \operatorname{Id} - C_1\Gamma^*\Gamma.$$

When $\widehat{G}$ and $\widehat{S}$ do not depend on $(t, y)$, one immediately gets symmetrizers for (2.8), taking $\widehat{S}$ as Fourier multiplier. When $\widehat{S}$ depends on $(t, y)$, one can try to take as a symmetrizer the pseudo-differential operator

$$S_\gamma(x) = \widehat{S}(t, y, x, D_t, D_y, \gamma)$$

or a para-differential version of it. The problem is then to convert the *symbolic* estimates (2.10) into estimates for the *operators* $S(x)$. This strategy has been shown to be successful for hyperbolic problems ([Kr], [CP]) and small viscosity parabolic-hyperbolic problems ([MZ], [GMWZ1], [GMWZ2]).

A further reduction to constant coefficients in $x$ can be performed when $\widehat{G}$ converges at an exponential rate to a limit when $x$ tend to infinity:

$$(2.11) \qquad |\widehat{G}(x,\zeta) - \widehat{G}(\infty,\zeta)| \leq Ce^{-\theta x}$$

for some $\theta > 0$ and $\zeta$ in a neighborhood of $\underline{\zeta}$. Then:

**Lemma 2.2 ([MZ]).** *There is a neighborhood $\omega$ of $\underline{\zeta}$ and there is a matrix $\mathcal{W}$ defined and $C^\infty$ on $[0,\infty[\times\omega$ such that*

  *i) $\mathcal{W}^{-1}$ is uniformly bounded and there is $\theta_1 > 0$ such that*

$$(2.12) \qquad |\mathcal{W}(x,\zeta) - \mathrm{Id}| \leq Ce^{-\theta_1 x})$$

  *ii) $\mathcal{W}$ satisfies*

$$(2.13) \qquad \partial_x \mathcal{W}(x,\zeta) = \widehat{G}(x,\zeta)\mathcal{W}(x,\zeta) - \mathcal{W}(x,\zeta)\mathcal{G}(\infty,\zeta) .$$

The substitution

$$v(x) = \mathcal{W}^{-1}(x,\zeta)\widehat{u}_\gamma(x), \qquad g(x) = \mathcal{W}^{-1}(x,\zeta)\widehat{f}_\gamma(x)$$

transforms (2.9) into

$$(2.14) \qquad \partial_x v = \widehat{G}(\infty,\zeta)v + g, \qquad \Gamma_1(\zeta)v(0) = 0,$$

with $\Gamma_1(\zeta) := \Gamma \mathcal{W}^{-1}(0,\zeta)$.

All this motivates the construction of symmetrizers for constant coefficients systems (1.1) with boundary conditions (1.2). Note that the reduction to constant coefficients above leads to frequency dependent boundary conditions. In the sequel, we concentrate on hyperbolic-parabolic systems arising in the study of boundary layers or shock profiles.

## 3 The main result

Consider a first order $N \times N$ system

$$(3.1) \qquad L_1(p,\partial)u := \partial_t + \sum_{j=1}^{d} A_j(p)\partial_j$$

The variables are $t \in \mathbb{R}$, $y \in \mathbb{R}^{d-1}$ and $x \in \mathbb{R}_+ := ]0, +\infty[$. Moreover, $\partial_j = \partial_{y_j}$ for $j \in \{1, \ldots, d-1\}$ and $\partial_d = \partial_x$. In addition, $p \in \mathbb{R}^M$ is a place holder for the variables themselves, the unknown in case of quasilinear systems and all the other parameters such as source terms or control.

Next, we consider a parabolic viscous perturbation of (3.1)

(3.2) $$L(p, \partial) = L_1(p, \partial) - L_2(p, \partial)$$

with

(3.3) $$L_2(p, \partial) := \sum_{1 \leq j,k \leq d} \partial_j \big(B_{j,k}(p) \partial_k u\big).$$

This kind of systems arises in the study of boundary layers or shock profiles, as linearized form of viscous perturbations of hyperbolic conservation laws, see e.g. [MZ] [GMWZ1], [GMWZ2]. In [MZ] [GMWZ2] we considered rather the small viscosity problem $L_1 - \varepsilon L_2$. Rescaling the variables as $t = \varepsilon \tilde{t}, y = \varepsilon \tilde{y}, x = \varepsilon \tilde{x}$ reduces to $\varepsilon = 1$ and the rescaled variables are much more adapted to the geometrical analysis of the symbol.

**Assumption 3.1.**

(H0) *The $A_j$ and $B_{j,k}$ are $C^\infty$ functions of $p$ on the open subset $\mathcal{O}$ of $\mathbb{R})M$, with values in the space of $N \times N$ real matrices,*

(H1) *For all $p \in \mathcal{O}$, the eigenvalues of $\sum \xi_j A_j(p)$ are real and semi-simple and have constant multiplicities for $p \in \mathcal{O}$ and $\xi \in \mathbb{R}^d \setminus \{0\}$.*

(H2) *There is $c > 0$ such that for all $p \in \mathcal{O}$ and $\xi \in \mathbb{R}^d$ the eigenvalues of $i\sum_{j=1}^d \xi_j A_j(u) + \sum_{j,k=1}^d \xi_j \xi_k B_{j,k}(u)$ satisfy $\operatorname{Re} \mu \geq c|\xi|^2$.*

(H3) *For all $p \in \mathcal{O}$, there holds $\det A_d(p) \neq 0$.*

We perform a Laplace-Fourier transform with respect to $t$ and a Fourier transform with respect to the variables $y$. This amounts to performing the substitution

(3.4) $$\partial_t \mapsto i\tau + \gamma, \quad \partial_{y_j} \mapsto i\eta_j$$

in the definition of $L$, with $\tau \in \mathbb{R}$, $\gamma \in \mathbb{R}$ and $\eta = (\eta_1, \ldots, \eta_{d-1}) \in \mathbb{R}^{d-1}$. In applications, to study the forward propagation, dictated by the sign of $L_2$, one restricts $\gamma$ to be nonnegative. After this Laplace-Fourier transform, $L$ becomes a second order differential operator in $x$ depending on the parameters $p$ and $\zeta = (\tau, \eta, \gamma) \in \mathbb{R}^{d+1}$. We reduce this system to first order introducing the variable $U = (u, \partial_x u)$ with values in $\mathbb{R}^{2N}$. This leads to the following $2N \times 2N$ system

$$(3.5) \qquad \mathcal{L} := \mathrm{Id}\partial_x - \mathcal{G}(p, \zeta), \qquad \mathcal{G}(p, \zeta) = \begin{pmatrix} 0 & \mathrm{Id} \\ \mathcal{M} & \mathcal{A} \end{pmatrix},$$

where

$$\begin{cases} \mathcal{A}(p, \zeta) = (B_{d,d})^{-1}\Big(A_d - \sum_{j=1}^{d-1} i\eta_j (B_{j,d} + B_{d,j})\Big) \\ \mathcal{M}(p, \zeta) = (B_{d,d})^{-1}\Big((i\tau + \gamma) + \sum_{j=1}^{d-1} i\eta_j A_j + \sum_{j,k=1}^{d-1} \eta_j \eta_k B_{j,k}\Big). \end{cases}$$

and the matrices $A_j$ and $B_{j,k}$ are evaluated at $p$.

The Assumption 3.1 immediately implies the following result (see [MZ] for a detailed proof).

**Lemma 3.2 (Spectral analysis of $\mathcal{G}$).**

*i) When $\zeta \neq 0$ and $\gamma \geq 0$, $\mathcal{G}(p, \zeta)$ has $N$ eigenvalues, counted with their multiplicities, in $\mathrm{Re}\,\mu > 0$ and $N$ eigenvalues in $\mathrm{Re}\,\mu < 0$.*

*ii) When $\zeta = 0$, $\mathcal{G}(p, 0)$ has $0$ as a semi-simple eigenvalue, of multiplicity $N$. The nonvanishing eigenvalues are those of $(B_{d,d})^{-1} A_d$.*

For $\zeta \neq 0$ and $\gamma \geq 0$ we denote by $\mathbb{E}_-(p, \zeta)$ [resp. $\mathbb{E}_+(p, \zeta)$] the spectral subspace of $\mathcal{G}(p, \zeta)$ associated to the spectrum in the left hand half space $\{\mathrm{Re}\,\mu < 0\}$ [resp. $\{\mathrm{Re}\,\mu > 0\}$]. Both have dimension $N$ and depend smoothly on $p$ and $\zeta$ as long as $\zeta \neq 0$ and $\gamma \geq 0$. To analyze the behavior of this spaces when $\zeta$ tends to zero, we introduce polar coordinates $\zeta = \rho\check\zeta$ with $\rho = |\zeta|$ and set

$$(3.6) \qquad \check{\mathbb{E}}(p, \check\zeta, \rho) := \mathbb{E}_-(p, \rho\check\zeta), \qquad \text{for } \rho > 0.$$

We can now state the main result of this paper:

**Theorem 3.3.** *Under Assumption 3.1 the linear bundle $\check{\mathbb{E}}(p, \check\zeta, \rho)$ has a continuous extension to $\rho = 0$, $\check\zeta \in \mathbb{R}^{d+1} \setminus \{0\}$ with $\check\gamma \geq 0$ and $p \in \mathcal{O}$.*

**Remark 3.4.** As already said, this result was already known under additional hypothesis [Rou], [ZS], [Z]. In [MZ] we used that the $\check{\mathbb{E}}_-(p, \check\gamma, \rho)$ is continuous up to $\rho = 0$ when $\check\gamma > 0$ and that $\check{\mathbb{E}}(p, \check\zeta, 0)$ has the form

$$(3.7) \qquad \check{\mathbb{E}}(p, \check\zeta, 0) = \mathbb{E}_-^{hyp}(p, \check\zeta) \oplus \mathbb{E}_-^{par}(p, 0),$$

see section 4 below. Moreover, $\mathbb{E}_-^{hyp}$ is the negative subspace corresponding to the hyperbolic system $L_1$. Following Kreiss [Kr] (see also [CP] and [Mé3] for the extension to hyperbolic systems of constant multiplicity) $\mathbb{E}_-^{hyp}(p, \check\zeta)$ can be extended continuously up to $\check\gamma = 0$. Therefore, in [MZ] we used the existence of the double limit

$$\lim_{\check\gamma \to 0} \lim_{\rho \to 0} \check{\mathbb{E}}_-(p, \check\zeta, \rho)$$

to establish desired linear decay estimates. The Theorem 3.3 asserts the existence of the limit as $\check\gamma$ and $\rho$ tend to zero, independently.

**Remark 3.5.** Consider a boundary condition $\Gamma(p, \zeta)U = 0$ for (3.5) where $\Gamma$ is a smooth function of $p$ and $\zeta$. The weak stability condition is given by (1.4) and the strong stability condition reads

(3.8) $$|D(p, \zeta)| = |\det(\mathbb{E}_-(p, \zeta), \ker \Gamma(p, \zeta))| \geq c$$

for some constant $c > 0$. When it holds, passing to the limit in $\rho$ implies that $|\det(\check{\mathbb{E}}_-(p, \check\zeta, 0), \ker \Gamma(p, 0))| \geq c$ for $\check\gamma > 0$. Next, passing to the limit in $\check\gamma$ shows that this condition holds up to $\check\gamma = 0$. Thus,

(3.9) $$\mathbb{C}^{2N} = \check{\mathbb{E}}_-(p, \check\zeta, 0) \oplus \ker \Gamma(p, 0)$$

when $\check\gamma = 0$. Theorem 3.3 when $\check\gamma = 0$ and the known smoothness when $\check\gamma > 0$, imply that conversely, if the transversality 3.9 holds at some point $(\underline{p}, \underline{\check\zeta})$ with $\underline{\check\zeta} \neq 0$ and $\underline{\check\gamma} \geq 0$, then the uniform stability condition 3.8 holds for $(p, \check\zeta, \rho)$ in a neighborhood of $(\underline{p}, \underline{\check\zeta}, 0)$ with $\check\gamma \geq 0$ and $\rho \geq 0$.

## 4 Strategy of the proof

By Lemma 3.2, for $\zeta$ small the eigenvalues of $\mathcal{G}$ split into a group of $N$ eigenvalues close to 0 and a group of $N$ eigenvalues away from the imaginary axis. More precisely, given a point $\underline{p}$, there is a $C^\infty$ invertible matrix $\mathcal{V}$ defined for $(p, \zeta)$ in neighborhood $\omega$ of $(\underline{p}, 0)$ such that $\mathcal{G}_1 := \mathcal{V}^{-1}\mathcal{G}\mathcal{V}$ has the block diagonal form

(4.1) $$\mathcal{G}_1(p, \zeta) = \begin{pmatrix} H(p, \zeta) & 0 \\ 0 & P(p, \zeta) \end{pmatrix}$$

with $H(p, 0) = 0$ and $P(p, 0) = (B_{d,d})^{-1}A_d$. The eigenvalues of $P(p, \zeta)$ satisfy $|\operatorname{Re}\mu| \geq c$ for some $c > 0$ and

(4.2) $$H = -(A_d)^{-1}\Big((i\tau + \gamma)\operatorname{Id} + \sum_{j=1}^{d-1} i\eta_j A_j\Big) + O(|\zeta|^2).$$

Therefore, for $(p, \zeta) \in \omega$, with $\zeta \neq 0$ and $\gamma \geq 0$:

$$\tag{4.3} \mathbb{E}_-(p, \zeta) = \mathcal{V}(p, \zeta)\big(\mathbb{E}_-^H(p, \zeta) \oplus \mathbb{E}_-^P(p, \zeta)\big)$$

where $\mathbb{E}_-^H$ [resp. $\mathbb{E}_-^P$] is the spectral space associated to eigenvalues of $H$ [resp. $P$] in $\{\operatorname{Re}\mu < 0\}$.

On one hand, because the spectrum of $P$ does not intersect the imaginary axis, the definition of the spaces $\mathbb{E}_-^P$ extends smoothly to $(p, \zeta) \in \omega$. On the other hand, the eigenvalues of $H$ are small, and we now investigate how they split in the complex complex domain when $\zeta$ is small. In polar coordinates, there holds

$$\tag{4.4} H(p, \rho\check{\zeta}) = \rho \check{H}(p, \check{\zeta}, \rho) = \rho\big(H_0(p, \check{\zeta}) + O(\rho)\big)$$

with

$$\tag{4.5} H_0(p, \check{\zeta}) = -(A_d)^{-1}\Big(i\check{\tau} + \check{\gamma})\operatorname{Id} + \sum_{j=1}^{d-1} i\check{\eta} A_j\Big)$$

Therefore

$$\tag{4.6} \mathbb{E}_-(p, \rho\check{\zeta}) = \mathcal{V}(p, \rho\check{\zeta})\big(\mathbb{E}_-^{\check{H}}(p, \check{\zeta}, \rho) \oplus \mathbb{E}_-^P(p, \rho\check{\zeta})\big)$$

where $\mathbb{E}_-^{\check{H}}$ is the negative space associated to $\check{H}$.

The hyperbolicity assumption (H1) implies that $H_0(p, \check{\zeta})$ has no eigenvalues on the imaginary axis when $\check{\gamma} \neq 0$. This remains true for $\check{H}(p, \check{\zeta}, \rho)$ for $\rho$ small enough (depending on $\check{\gamma}$). In particular, for all $\check{\zeta} \in \mathbb{R}^{d+1}$ with $\check{\gamma} > 0$, the definition of the bundle $\mathbb{E}^{\check{H}}(p, \check{\zeta}, \rho)$ extends smoothly to a neighborhood of $(\underline{p}, \underline{\check{\zeta}}, 0)$ in $\mathbb{R}^M \times \mathbb{R}^{d+1} \times \mathbb{R}$. This remains true for $\mathbb{E}(p, \rho\check{\zeta})$ since $\mathbb{E}_-^P(p, \rho\check{\zeta})$ and $\mathcal{V}(p, \rho\check{\zeta})$ are smooth in $\rho\check{\zeta}$.

Thus, it remains to prove Theorem 3.3 on a neighborhood of a point $(\underline{p}, \underline{\check{\zeta}}, 0)$ where $\check{\gamma} = 0$ and $\check{\zeta} \neq 0$.

The strategy of the proof is as follows. We first construct explicitly spaces $\underline{\mathbb{E}}_-$ and $\underline{\mathbb{E}}_+$ such that

$$\tag{4.7} \mathbb{C}^{2N} = \underline{\mathbb{E}}_- \oplus \underline{\mathbb{E}}_+, \quad \dim \underline{\mathbb{E}}_- = \dim \underline{\mathbb{E}}_+ = N.$$

Indeed, extending Kreiss' analysis to hyperbolic systems of constant multiplicity as in [Mé3], one can show that the bundle $\mathbb{E}_-^{H_0}(p, \check{\zeta})$ extends continuously to $\check{\gamma} = 0$ and we take $\underline{\mathbb{E}}_-$ to be equal to the value of this extension at $(\underline{p}, \underline{\check{\zeta}})$, but we will not use this result, see below. Next the main point is to construct symmetrizers. In the next statement we denote by $\underline{\Pi}_\pm$ the projectors on $\underline{\mathbb{E}}_\pm$ associated to the decomposition (4.7).

**Theorem 4.1.** *For all $(\underline{p}, \underline{\check\zeta})$ with $\underline{\check\zeta} \ne 0$ and $\underline{\check\gamma} = 0$, there are spaces $\underline{\mathbb{E}}_\pm$ satisfying (4.7) and such that for all $\kappa \geq 1$ there are a neighborhood $\Omega$ of $(\underline{p}, \underline{\check\zeta})$ in $\mathbb{R}^M \times \mathbb{R}^{d+1} \times \mathbb{R}$, a $C^\infty$ mapping $\mathcal{S}$ from $\Omega$ to the space of $2N \times 2N$ matrices and a constant $c > 0$ such that for all $(p, \check\zeta, \rho) \in \Omega$,*

(4.8) $$\mathcal{S}(p, \check\zeta, \rho) = \mathcal{S}^*(p, \check\zeta, \rho),$$

(4.9) $$\big(\mathcal{S}(p, \check\zeta, \rho) U, U\big) \geq \kappa^2 |\underline{\mathrm{II}}_+ U|^2 - |\underline{\mathrm{II}}_- U|^2,$$

*and for all $(p, \check\zeta, \rho) \in \Omega$ with $\rho \geq 0$ and $\check\gamma \geq 0$:*

(4.10) $$\mathrm{Re}\,\big(\mathcal{S}(p, \check\zeta, \rho) \mathcal{G}(p, \rho\check\zeta) U, U\big) \geq c\rho(\check\gamma + \rho)|U|^2.$$

*Proof of Theorem 3.3 assuming Theorem 4.1.*

**a)** Consider $\kappa > 2$ and $\Omega$ given by Theorem (4.1). For $(p, \check\zeta, \rho) \in \Omega$ with $\check\gamma \geq 0$ and $\rho > 0$, let $U \in \mathbb{E}_-(p, \rho\check\zeta)$. Then

$$\widetilde{U}(x) = e^{x\mathcal{G}(p, \rho\check\zeta)} U$$

is exponentially decaying at $+\infty$ and satisfies $\partial_x \widetilde{U} = \mathcal{G}(p, \rho\check\zeta)\widetilde{U}$. Therefore multiplying by $\mathcal{S}(p, \check\zeta, \rho)$ and integrating by parts as in Lemma 2.1, yields thanks to (4.8):

$$\big(\mathcal{S}U, U\big) + 2\mathrm{Re} \int_0^\infty \big(\mathcal{S}\mathcal{G}\widetilde{U}(x), \widetilde{U}(x)\big) dx = 0$$

By (4.10), the integral is nonnegative. Therefore $(\mathcal{S}U, U)$ is nonpositive which by (4.9) implies that $\kappa|\underline{\mathrm{II}}_+ U| \leq |\underline{\mathrm{II}}_- U|$. Thus

(4.11) $$\forall U \in \mathbb{E}_-(p, \rho\check\zeta) : \qquad |\underline{\mathrm{II}}_+ U| \leq \frac{1}{\kappa - 1}|U|.$$

This implies that the mapping $\underline{\mathrm{II}}_-$ from $\mathbb{E}_-(p, \rho\check\zeta)$ into $\underline{\mathbb{E}}_-$ is one to one and since both spaces have dimension $N$, it is a bijection. Therefore, there is a mapping $\Phi(p, \check\zeta, \rho)$ from $\underline{\mathbb{E}}_-$ to $\underline{\mathbb{E}}_+$ such that

(4.12) $$\mathbb{E}_-(p, \rho\check\zeta) = \big\{ u + A(p, \check\zeta, \rho)u \,:\, u \in \underline{\mathbb{E}}_- \big\}$$

and

(4.13) $$\forall u \in \underline{\mathbb{E}}_- : \qquad |A(p, \check\zeta, \rho) u| \leq \frac{1}{\kappa - 2}|u|.$$

Since $\kappa$ is arbitrarily large, this proves that

(4.14) $$\underline{\mathbb{E}}_- = \widetilde{\lim}\, \mathbb{E}_-(p, \rho\check\zeta),$$

where $\widetilde{\lim}$ means that $(p, \check\zeta, \rho)$ tends to $(\underline{p}, \underline{\check\zeta}, 0)$ with $\rho > 0$, $\check\gamma \geq 0$.

**b)** Since this is true for all $(p, \check{\zeta})$ this proves that $\widetilde{\mathbb{E}}_-(p, \check{\zeta}, \rho) := \mathbb{E}_-(p, \rho\check{\zeta})$ extends to points $(p, \check{\zeta}, 0)$ with $\check{\gamma} \geq 0$. We denote by $\widetilde{\mathbb{E}}_-(p, \check{\zeta}, \rho)$ this extension and (4.14) implies that for all $(\underline{p}, \underline{\check{\zeta}})$ with $\underline{\check{\gamma}} \geq 0$, there holds

(4.15) $$\widetilde{\mathbb{E}}_-(\underline{p}, \underline{\check{\zeta}}, 0) = \widetilde{\lim} \, \widetilde{\mathbb{E}}_-(p, \check{\zeta}, \rho).$$

Consider again a given point $(\underline{p}, \underline{\check{\zeta}})$ with $\underline{\check{\gamma}} = 0$. For $\kappa > 2$, let $\Omega$ be given by Theorem 4.1. For all $(p', \check{\zeta}', 0) \in \Omega$, thanks to (4.15), we can let $(p, \check{\zeta}, \rho)$ tend $(p', \check{\zeta}', 0)$ in the sense of $\widetilde{\lim}$. Therefore, (4.11) implies that for all $(p, \check{\zeta}, \rho) \in \Omega$ with $\check{\gamma} \geq 0$ and $\rho \geq 0$:

$$\forall U \in \widetilde{\mathbb{E}}_-(p, \check{\zeta}, \rho) : \qquad |\underline{\mathrm{II}}_+ U| \leq \frac{1}{\kappa - 1} |U|.$$

Arguing as before, this implies that

$$\widetilde{\mathbb{E}}_-(\underline{p}, \underline{\check{\zeta}}, 0) = \underline{\mathbb{E}}_- = \lim \widetilde{\mathbb{E}}_-(p, \check{\zeta}, \rho),$$

where the limit is taken for $(p, \check{\zeta}, \rho)$ tending $(\underline{p}, \underline{\check{\zeta}}, 0)$ with $\rho \geq 0$, $\check{\gamma} \geq 0$.

This means that the bundle $\widetilde{\mathbb{E}}_-(p, \check{\zeta}, \rho)$ is continuous in $(p, \check{\zeta}, \rho)$ for $\check{\gamma} \geq 0$ and $\rho \geq 0$. $\square$

## 5 Proof of Theorem 4.1

Though not stated explicitly, the result of Theorem 4.1 is given in Appendix A, section 3, of [MZ]. For the convenience of the reader, we give here the main steps of the proof.

**a)** We first remark that it is sufficient to prove the theorem for $\kappa$ large enough. Next, we can replace (4.9) by the weaker condition

(5.1) $$\left(\underline{\mathcal{S}} U, U\right) \geq \kappa^2 |\underline{\mathrm{II}}_+ U|^2 - |\underline{\mathrm{II}}_- U|^2,$$

for $\underline{\mathcal{S}} = \mathcal{S}(\underline{p}, \underline{\check{\zeta}}, 0)$. Indeed, if (5.1) holds, then, by continuity of $\mathcal{S}$, on a possibly smaller neighborhood of $(\underline{p}, \underline{\check{\zeta}}, 0)$, there holds

$$\left(\mathcal{S}(p, \check{\zeta}, \rho) U, U\right) \geq \frac{1}{2} \kappa^2 |\underline{\mathrm{II}}_+ U|^2 - 2 |\underline{\mathrm{II}}_- U|^2,$$

and therefore $\frac{1}{2} \mathcal{S}$ satisfies (4.8) for $\kappa/2$ and also (4.8) and (4.10) with the constant $c/2$.

**b)** Recall that $\mathcal{G}_1 = \mathcal{V}^{-1}\mathcal{G}\mathcal{V}$ has the form (4.1). We construct spaces $\underline{\mathbb{E}}^H_\pm$ and $\underline{\mathbb{E}}^P_\pm$ such that

$$(5.2) \quad \mathbb{C}^N = \underline{\mathbb{E}}^H_- \oplus \underline{\mathbb{E}}^H_+, \quad \mathbb{C}^N = \underline{\mathbb{E}}^P_- \oplus \underline{\mathbb{E}}^P_+ \quad \dim \underline{\mathbb{E}}^H_- + \dim \underline{\mathbb{E}}^P_+ = N.$$

Denote by $\underline{\Pi}^H_\pm$ and $\underline{\Pi}^P_\pm$ the corresponding projectors. Suppose that for all $\kappa$ large enough there are symmetrizers $S^H$ and $S^P$ defined on some neighborhood of $(p,\check{\zeta})$ and such that

$$(5.3) \qquad S^H = (S^H)^*$$
$$(5.4) \qquad (\underline{S}^H U, U)) \geq \kappa^2 |\underline{\Pi}^H_+ U|^2 - |\underline{\Pi}^H_- U|^2,$$
$$(5.5) \qquad \operatorname{Re} S^H \check{H} \geq c(\check{\gamma} + \rho)\operatorname{Id}, \quad \text{for } \rho \geq 0,\ \check{\gamma} \geq 0,$$

and

$$(5.6) \qquad S^P = (S^P)^*$$
$$(5.7) \qquad (\underline{S}^P U, U)) \geq \kappa^2 |\underline{\Pi}^P_+ U|^2 - |\underline{\Pi}^P_- U|^2,$$
$$(5.8) \qquad \operatorname{Re} S^P P \geq c\operatorname{Id}.$$

Here, we have denoted by $\underline{S}^H$ and $\underline{S}^P$ the value of $S^H$ and $S^P$ respectively at the base point $(p,\check{\zeta},0)$.

Then we choose $\underline{\mathbb{E}}_\pm = \mathcal{V}(\underline{\mathbb{E}}^H_\pm \oplus \underline{\mathbb{E}}^P_\pm)$ and $\mathcal{S} = (\mathcal{V}^{-1})^* \mathcal{S}_1 \mathcal{V}^{-1}$ with

$$(5.9) \qquad \mathcal{S}_1(p,\zeta) = \delta \begin{pmatrix} S^H & 0 \\ 0 & S^P \end{pmatrix}$$

By (5.3) and (5.6), $\mathcal{S}$ is clearly self adjoint and (5.5) and (5.8) imply (4.10). In addition, for $U = \mathcal{V}(U^H \oplus U^P)$ one has

$$(\underline{\mathcal{S}}U, U) = \delta(S^H U^H, U^H) + \delta(S^P U^P, U^P).$$

Because $\underline{\Pi}_\pm = \mathcal{V}(\underline{\Pi}^H_\pm \oplus \underline{\Pi}^P_\pm)\mathcal{V}^{-1}$, there is $C$ such that

$$|\underline{\Pi}^H_+ U^H|^2 + |\underline{\Pi}^P_+ U^P|^2 \geq \frac{1}{C}|\underline{\Pi}_+ U|^2$$
$$|\underline{\Pi}_- U|^2 \leq C\left(|\underline{\Pi}^H_- U^H|^2 + |\underline{\Pi}^P_- U^P|^2\right)$$

Therefore, choosing $\delta = 1/C$, (5.4) and (5.7) imply (5.1) with $\kappa$ replaced by $\kappa/C$ and it is sufficient to construct the symmetrizers $S^H$ and $S^P$.

**c)** The construction of $S^P$ is easy. Because the eigenvalues of $P$ do not belong to the imaginary axis, one can perform a further reduction (or choose $\mathcal{V}$ in (4.1)) such that for $(p, \zeta)$ in a neighborhood of $(\underline{p}, 0)$:

$$P(p, \zeta) = \begin{pmatrix} P_+(p, \zeta) & 0 \\ 0 & P_-(p, \zeta) \end{pmatrix}$$

where the spectrum of $P_\pm$ is contained in $\{\pm \operatorname{Re} \mu > 0\}$. In this basis, the spaces $\mathbb{E}_\pm(p, \zeta)$ are constant. we call them $\underline{\mathbb{E}}_\pm$. They are generated by the first and last vectors of the basis corresponding to the blocks $P_+$ and $P_-$ respectively.

There are self-adjoint matrices $S_\pm^P$ and $c > 0$ such that

$$\pm \operatorname{Re} S_\pm^P P_\pm \geq c \operatorname{Id}, \quad S_+^P \geq \operatorname{Id}, \quad S_-^P \leq \operatorname{Id}.$$

For instance, one can choose with appropriate positive constants $C_\pm$

$$S_\pm^P = C_\pm \int_0^\infty \left(e^{\mp t P_\pm}\right)^* e^{\mp t P_\pm}\, dt.$$

Then the symmetrizers

$$S^P(p, \zeta) = \begin{pmatrix} \kappa^2 S_+^P(p, \zeta) & 0 \\ 0 & -S_-^P(p, \zeta) \end{pmatrix}$$

satisfy (5.6) to (5.8).

**Remark 5.1.** The symmetrizers $S^P$ constructed above are smooth functions of $(p, \zeta)$. This was important in [MZ] in the quantization of the symbolic calculus.

**d)** We now turn to the construction of $S^H$. As already mentioned in (4.4), $\check{H}$ is a perturbation of $H_0$. We denote by $\underline{\mu}_k = i\underline{\xi}_k$ the distinct eigenvalues of $\underline{H}_0 := H_0(\underline{p}, \check{\underline{\zeta}})$ located on the imaginary axis. There is $\delta > 0$, a neighborhood $\Omega_0$ of $(\underline{p}, \check{\underline{\zeta}}, 0)$ and a smooth matrix $V(p, \check{\zeta}, \rho)$ such that

(5.10) $$V^{-1} \mathcal{H} V = \begin{bmatrix} \mathcal{H}_1 & \cdots & 0 \\ \vdots & \ddots & \vdots \\ 0 & \cdots & \mathcal{H}_k \end{bmatrix}$$

such that each block $\mathcal{H}_k$ has its spectrum either in $\{|\operatorname{Re} \mu| \geq 2\delta\}$ or in the ball of radius $\delta$ centered at $\underline{\mu}_k$. Moreover, we can assume that the balls of radius $2\delta$ centered at the $\mu_k$ do not intersect each other.

We denote by $N_k$ the dimension of the block $\mathcal{H}_k$. By Assumption 3.1, the $\mathcal{H}_k$ have no eigenvalues on the imaginary axis when $\rho \geq 0$, $\check\gamma \geq 0$ and $\rho + \check\gamma > 0$. Therefore, the number of eigenvalues of $\mathcal{H}_k$ in $\{\pm \operatorname{Re}\mu > 0\}$ is constant for $\rho \geq 0$, $\check\gamma \geq 0$ and $\rho + \check\gamma > 0$. We denote it by $N_{k,\pm}$. Because the total number of eigenvalues of $\mathcal{G}$ in $\{\operatorname{Re}\mu < 0\}$ is $N$, we have

$$\sum N_{k,-} + \dim \mathbb{E}_-^P = N. \tag{5.11}$$

Arguing as in step b), it is sufficient to construct spaces $\underline{\mathbb{E}}_\pm^k$ satisfying

$$\mathbb{C}^{N_k} = \underline{\mathbb{E}}_+^k \oplus \underline{\mathbb{E}}_-^k, \quad \dim \underline{\mathbb{E}}_-^k = N_{k,-}, \tag{5.12}$$

and for all $\kappa$ large enough, symmetrizers $\mathcal{S}^k$, $C^\infty$ on a neighborhood of $(\underline{p}, \underline{\check\zeta}, 0)$, such that

$$\mathcal{S}^k = (\mathcal{S}^k)^* \tag{5.13}$$

$$(\underline{\mathcal{S}}^k U, U)) \geq \kappa^2 |\underline{\Pi}_+^k U|^2 - |\underline{\Pi}_-^k U|^2, \tag{5.14}$$

$$\operatorname{Re} \mathcal{S}^k \mathcal{H}^k \geq c(\check\gamma + \rho)\operatorname{Id}, \quad \text{for } \rho \geq 0, \ \check\gamma \geq 0. \tag{5.15}$$

In (5.14), $\underline{\mathcal{S}}^k = \mathcal{S}_k(\underline{p}, \underline{\check\zeta}, 0)$ and $\underline{\Pi}_\pm^k$ are the projectors associated to the decomposition (5.12).

**e)** If the spectrum of $\mathcal{H}^k$ lies in $\{\operatorname{Re}\mu > 2\delta\}$ [resp. $\{\operatorname{Re}\mu < -2\delta\}$], then arguing as in step c), we set $\underline{\mathbb{E}}_-^k = \{0\}$ [resp. $\underline{\mathbb{E}}_-^k = \mathbb{C}^{N_k}$] and $\mathcal{S}^k = \kappa^2 S^k$ [resp. $\mathcal{S}^k = -S^k$] where $S^k(p, \check\zeta, \rho)$ is a self adjoint matrix, defined and $C^\infty$ on a neighborhood of $(\underline{p}, \underline{\check\zeta}, 0)$ and such that $\operatorname{Re} S^k \geq \operatorname{Id}$ and $\operatorname{Re} S^k \mathcal{H}^k$ is positive definite [resp. $\operatorname{Re} S^k \leq \operatorname{Id}$ and $\operatorname{Re} S^k \mathcal{H}^k$ is negative definite].

**f)** Consider now the case where the spectrum of $\mathcal{H}^k$ is contained in the ball of radius $\delta$ centered at $\underline{\mu}_k = i\underline{\xi}_k$. We note that $i\underline{\xi}_k$ is an eigenvalue of $\underline{H}_0$ if an only if $-\underline{\check\tau} = \lambda(\underline{p}, \underline{\check\eta}, \underline{\xi}_k)$ where $\lambda(p, \eta, \xi)$ is one among the eigenvalues of $\sum \eta_j A_j(p) + \xi A_d(p)$.

Suppose first that $\underline{\xi}_k$ is a simple root of the eigenvalue equation, that is

$$\partial_\xi \lambda(\underline{p}, \underline{\check\eta}, \underline{\xi}_k) \neq 0. \tag{5.16}$$

Denote by $\alpha_k$ the multiplicity of $\lambda$ (see Assumption (H1)). Then, according to [Mé3], $\mathcal{H}^k(p, \check\zeta, 0)$ is a scalar matrix $q^k(p, \check\zeta)\operatorname{Id}$ and by Lemma 2.10 of [MZ]

$$\mathcal{H}^k(p, \check\zeta, \rho) = q^k(p, \check\zeta)\operatorname{Id} + \rho \mathcal{R}^k(p, \check\zeta, \rho). \tag{5.17}$$

Moreover, $q_k$ is purely imaginary when $\check\gamma = 0$, $\dot q_k := \partial_{\check\gamma} \operatorname{Re} q^k(\underline p, \underline{\check\zeta})$ does not vanish and $\dot q_k \operatorname{Re} \mathcal R^k(\underline p, \underline{\check\zeta}, 0)$ is definite positive.

When $\dot q_k > 0$ [resp. $\dot q_k < 0$], for $\check\gamma > 0$ small, the eigenvalue $q_k$ has a positive [resp. negative] real part. Thus we set $\underline{\mathbb E}^k_- = \{0\}$ [resp. $\underline{\mathbb E}^k_- = \mathbb C^{N_k}$]. Next we choose $\mathcal S^k = \kappa^2 S^k$ [resp. $\mathcal S^k = -S^k$] where $S^k(p, \check\zeta, \rho)$ is a self adjoint matrix, defined and $C^\infty$ on a neighborhood of $(\underline p, \underline{\check\zeta}, 0)$ and such that $\operatorname{Re} S^k \ge \operatorname{Id}$ and $\operatorname{Re} S^k \mathcal R^k$ is positive definite [resp. $\operatorname{Re} S^k \le \operatorname{Id}$ and $\operatorname{Re} S^k \mathcal R^k$ is negative definite].

**g)** We now come to the most difficult part of the construction, when $\underline\xi_k$ is a multiple root of the eigenvalue equation, that is, when there is $\nu_k \ge 2$ such that
$$\partial_\xi \lambda = \ldots \partial_\xi^{\nu_k - 1} \lambda = 0, \qquad \partial_\xi^{\nu_k} \ne 0$$
at $(\underline p, \underline{\check\eta}, \underline\xi_k)$. (Note that, since $\lambda$ is real analytic in $\xi$, there is always such an integer $\nu_k$). The case $\nu_k = 1$ is (5.16)).

¿From [Kr] in the strictly hyperbolic case and from [Mé3] for the extension to constant multiplicity, we know that there is a smooth matrix $\mathcal V^k_0(p, \check\zeta)$ such that

$$(5.18) \quad \mathcal Q^k(p, \check\zeta) := (\mathcal V^k_0)^{-1}(p, \check\zeta) \mathcal H^k(p, \check\zeta, 0) \mathcal V^k_0(p, \check\zeta) = \begin{bmatrix} Q_k & \cdots & 0 \\ \vdots & \ddots & \vdots \\ 0 & \cdots & Q_k \end{bmatrix},$$

where the subblock $Q_k$ is $\nu_k \times \nu_k$ matrices and there are $\alpha_k$ such blocks, where $\alpha_k$ is the multiplicity of the eigenvalue $\lambda$. In particular the dimension of $\mathcal H^k$ is $N_k = \alpha_k \nu_k$. Moreover,

$$(5.19) \qquad \underline Q_k := Q_k(\underline p, \underline{\check\zeta}) = i \begin{bmatrix} \underline\xi_k & 1 & 0 \\ 0 & \underline\xi_k & \ddots & 0 \\ & \ddots & \ddots & 1 \\ & & \cdots & \underline\xi_k \end{bmatrix}.$$

In addition, with Ralston's lemma [Ral], we can assume that only the first column of $Q_k$ does not vanish:

$$(5.20) \qquad Q_k(\underline p, \underline{\check\zeta}, 0) = \begin{bmatrix} * & 0\ldots 0 \\ \vdots & 0\ldots 0 \\ q_k & 0\ldots 0 \end{bmatrix},$$

where $q_k$ denotes the lower left hand corner of $Q_k$, that $Q_k$ has purely imaginary coefficients when $\check\gamma = 0$ and that $\dot q_k := \partial_{\check\gamma} \operatorname{Re} q_k(\underline p, \underline{\check\zeta}) \ne 0$.

Next from Lemma 2.10 of [MZ], we know that there is a smooth matrix $\mathcal{V}^k$ on a neighborhood of $(\underline{p}, \check{\zeta}, 0)$, which extends $\mathcal{V}_0$, such that

$$(5.21) \qquad \mathcal{H}_1^k := (\mathcal{V}^k)^{-1} \mathcal{H}^k \mathcal{V}^k = \mathcal{Q}^k(p, \check{\zeta}) + \rho \mathcal{R}^k(p, \check{\zeta}, \rho),$$

where the matrix $\mathcal{R}^k$ has the following decomposition in $\alpha_k \times \alpha_k$ blocks $R_{p,q}^k$ of size $\nu_k \times \nu_k$:

$$(5.22) \qquad \mathcal{R}_k = \begin{bmatrix} R_{1,1}^k & \cdots & R_{1,\alpha_k}^k \\ \vdots & \ddots & \vdots \\ R_{\alpha_k,1}^k & \cdots & R_{\alpha_k,\alpha_k}^k \end{bmatrix}.$$

Moreover, only the first column of the $R_{p,q}^k$ does not vanish:

$$(5.23) \qquad R_{p,q}^k(\underline{p}, \check{\zeta}, 0) = \begin{bmatrix} * & 0 \ldots 0 \\ \vdots & 0 \ldots 0 \\ r_{p,q}^k & 0 \ldots 0 \end{bmatrix},$$

where $r_{p,q}^k$ denotes the lower left hand entry of $R_{p,q}^k$. In addition, denoting by $R_k^{\flat}$ the $\alpha_k \times \alpha_k$ matrix with entries $r_{p,q}^k$, the matrix $\dot{q}_k \mathrm{Re}\, R_k^{\flat}(\underline{p}, \check{\zeta}, 0)$ is definite positive.

Arguing as in step a) and b), we are reduced to construct spaces and symmetrizers for each block $\mathcal{H}_1^k$.

**h)** From [Kr] [CP], or simply from a direct analysis of the model case

$$(5.24) \qquad i \begin{bmatrix} \underline{\xi}_k & 1 & & 0 \\ 0 & \underline{\xi}_k & \ddots & \\ & & \ddots & 1 \\ & \cdots & & \underline{\xi}_k \end{bmatrix} + \check{\gamma} \begin{bmatrix} * & 0 & \cdots & 0 \\ * & 0 & \cdots & 0 \\ \vdots & 0 & \cdots & 0 \\ \dot{q}_k & 0 & \cdots & 0 \end{bmatrix},$$

we know that the candidate for the limit negative space for this block is the space generated by the first $\beta_k$ vectors of the canonical basis:

$$(5.25) \qquad \mathbf{E}_-^k(\underline{p}, \check{\zeta}) = \mathbb{C}^{\beta_k} \times \{0\}^{\nu_k - \beta_k}$$

where

$$(5.26) \qquad \beta_k = \begin{cases} \frac{1}{2}\nu_k & \text{when } \nu_k \text{ is even}, \\ \frac{1}{2}(\nu_k \pm 1) & \text{when } \nu_k \text{ is odd} \end{cases} \quad \text{and } \mp \dot{q}_k > 0.$$

We also introduce
$$\mathbf{E}_+^k(\underline{p},\underline{\check{\zeta}}) = \{0\}^{\beta_k} \times \mathbb{C}^{\nu_k-\beta_k}$$
Thus
$$\mathbb{C}^{\nu_k} = \mathbf{E}_-^k \oplus \mathbf{E}_+^k .$$
In the block decomposition of $\mathbb{C}^{N_k}$ into $\alpha_k$ factors $\mathbb{C}^{\nu_k}$, let
$$\underline{\mathbb{E}}_\pm^k(p,\check{\zeta},0) = \mathbf{E}_\pm^k(p,\check{\zeta}) \oplus \cdots \oplus \mathbf{E}_\pm^k(p,\check{\zeta}).$$

Note that $\beta_k$ is exactly the number of eigenvalues in $\{\operatorname{Re}\mu < 0\}$ of the model (5.24) for $\check{\gamma} > 0$. This is still true for $Q_k(p,\underline{\check{\tau}},\underline{\check{\eta}},\check{\gamma})$. Thus the number of eigenvalues in $\{\operatorname{Re}\mu < 0\}$ of $\mathcal{Q}^k$ is $\nu_k\beta_k$ when $\check{\gamma} >$ and therefore the number of eigenvalues of $\mathcal{H}_1^k$ in $\{\operatorname{Re}\mu < 0\}$ when $\rho \geq 0$, $\check{\gamma} \geq 0$ and $\rho + \check{\gamma} > 0$, which is constant on a neighborhood of $(\underline{p},\underline{\check{\zeta}},0)$, is
$$N_k = \nu_k\beta_k = \dim \underline{\mathbb{E}}_-^k .$$
Therefore, the condition (5.12) is satisfied.

We now proceed to the construction of the symmetrizers. We construct $\mathcal{S}^k$ in the block decomposition of $\mathcal{H}_1^k$

(5.27) $$\mathcal{S}_k = \delta \begin{bmatrix} S_k & 0 & \\ 0 & S_k & \\ & & \ddots \end{bmatrix},$$

with

(5.28) $$S_k(p,\check{\zeta},\rho) = E_k + \widetilde{E}_k(p,\check{\zeta}) - i\gamma F_k - i\rho F_k',$$

where $E_k$ and $\widetilde{E}_k$ are real symmetric matrices, and $F_k$ and $H_k$ are real and skew symmetric. $E_k$ is constant and has the special form

$$E_k = \begin{bmatrix} 0 & \cdots & \cdots & 0 & e_{k,1} \\ \vdots & & & \iddots & e_{k,2} \\ \vdots & & \iddots & \iddots & \\ 0 & \iddots & \iddots & & \\ e_{k,1} & e_{k,2} & & & e_{k,\nu_k} \end{bmatrix}.$$

Moreover $\widetilde{E}_k(\underline{p},\underline{\check{\zeta}}) = 0$ and $F_k$ and $F_k'$ are constant.

The order of the construction is as follows. One first chooses $E_k$, $\widetilde{E}_k$ and $F_k$ as in [Kr] (see also [CP]) to construct symmetrizers for $Q_k$, that is for $\rho = 0$. The new part lies in the choice of $F_k'$.

1. Choice of $E_k$. First choose the real coefficient $e_{k,1}$ such that

(5.29) $$e_{k,1}\dot{q}_k \geq 3.$$

Next, the coefficients $e_{k,l}$ for $l > 1$ are chosen successively to achieve that there is $c > 0$ such that

(5.30) $$(E_k w, w) \geq c\Big(\kappa^2 |\underline{\pi}_+^k w|^2 - |\underline{\pi}_-^k w|^2\Big),$$

where $\underline{\pi}_\pm^k$ is the projection onto $\mathbf{E}_\pm^k$ in the decomposition $\mathbb{C}^{\nu_k} = \mathbf{E}_+^k \oplus \mathbf{E}_-^k$. (cf Lemma 5.6 in [CP], Chap. 7).

Next, following [CP] (cf (equation (5.5.3) in Chap. 7), (5.29) implies that there is a constant $C$ and a neighborhood of $(\underline{p}, \check{\underline{\zeta}})$ such that

(5.31) $$\mathrm{Re}\,\big(E_k \partial_\gamma Q_k(\underline{p}, \check{\underline{\zeta}}) w, w\big) \geq 2|w_1|^2 - C|w'|^2,$$

with $w_1$ the first component of $w \in \mathbb{C}^{\nu_k}$ and $w' \in \mathbb{C}^{\nu_k - 1}$ denotes the other components.

2. Choice of $\widetilde{E}_k$. Recall that $\underline{Q}_k = i(\underline{\xi}_k \mathrm{Id} + N_k)$ where $N_k$ is the Jordan matrix of size $\nu_k$ (see (5.19)). The form of $E_k$ is chosen so that $E_k(\underline{\xi}_k \mathrm{Id} + N_k)$ is real and symmetric. Next, the real matrix $\widetilde{E}_k(\underline{p}, \check{\underline{\zeta}})$ is chosen so that such that $(E_k + \widetilde{E}_k)(\frac{1}{i}Q_k)$ is real and symmetric when $\check{\gamma} = 0$. This is achieved in [Kr] [CP] using the implicit function theorem and the property that $\frac{1}{i}Q_k$ is real when $\check{\gamma} = 0$.

3. Choice of $F_k$. Following [Kr] [CP], there is $F_k$ real and skew symmetric such that

$$\mathrm{Re}\,(F_k N_k w, w) \geq -|w_1|^2 + (C+1)|w'|^2.$$

where $C$ is the constant in (5.31). As a consequence, we have

(5.32) $$\begin{aligned}\mathrm{Re}\,\big((E_k + \widetilde{E}_k - i\check{\gamma} F_k)Q_k\big) &= \check{\gamma} D_k\,,\\ D_k(\underline{p}, \check{\underline{\zeta}}) &= \mathrm{Re}\,(E_k \partial_\gamma Q_k(\underline{p}, \check{\underline{\zeta}})) + \mathrm{Re}\,(F_k N_k) \geq \mathrm{Id}\,.\end{aligned}$$

4. Choice of $F_k'$. This is the new part detailed in the Appendix A of [MZ]. Denote by $\mathcal{E}_k$ the block diagonal matrix $\mathrm{Diag}(E_k)$. A vector $w \in \mathbb{C}^{N_k}$, $N_k = \nu_k \alpha_k$, is broken into $\alpha_k$ blocks $w_p \in \mathbb{C}^{\nu_k}$, with components denoted by $w_{p,a}$. We denote by $R_{p,q}^k$ the $\nu_k \times \nu_k$ blocks of $\mathcal{R}_k$ and by $R_{p,a,q,b}$ their entries. The entries of $E_k$ are denoted by $E_{a,b}$. Since $R_{p,a,q,b} = 0$ when $b > 1$, the special form of $E_k$ implies that

$$\begin{aligned}\mathrm{Re}\,(\mathcal{E}_k \mathcal{R}_k w, w) &= \mathrm{Re}\,\sum E_{a,c} R_{p,a,q,1} w_{q,1} \overline{w}_{p,c}\\ &= \mathrm{Re}\,\sum e_{k,1} r_{p,q} w_{q,1} \overline{w}_{p,1} + O(|w_{*,1}|\,|w_*'|)\end{aligned}$$

where $w_{*,1} \in \mathbb{C}^{\alpha_k}$ is the collection of the first components $w_{p,1}$, $w'_*$ the remainder components and $r_{p,q} = R_{p,\nu_k,q,1}$ the lower left hand corner entry of $R_{p,q}$. The matrix $\text{Re}\,(R_k^\flat)$ is definite, positive or negative according to the sign of $\dot{q}_k$, which is the sign of $e_{k,1}$ by (5.29). Thus, multiplying $E_k$ by some positive constant, we can achieve that in addition to (5.30) (5.32), the following inequality holds:

$$\text{(5.33)} \qquad \text{Re}\,\big(\mathcal{E}_k \mathcal{R}_k(\underline{p}, \check{\underline{\zeta}}, 0) w, w\big) \geq 2|w_{*,1}|^2 - C'|w'_*|^2.$$

for some constant $C' > 0$.

As in 3, there is $F'_k$ real and skew symmetric such that for all $w \in \mathbb{C}^{\nu_k}$:

$$\text{Re}\,(F'_k N_k w, w) \geq -|w_1|^2 + (C'+1)|w'|^2.$$

Thus, with $\mathcal{F}'_k = \text{Diag}(F'_k)$, $\mathcal{N}_k = \text{Diag}(N_k)$ and $w \in \mathbb{C}^{N_k}$:

$$\text{Re}\,(\mathcal{F}'_k \mathcal{N}_k w, w) \geq -|w_{*,1}|^2 + (C'+1)|w'_*|^2.$$

Therefore, with (5.33), we have

$$\text{(5.34)} \qquad \text{Re}\,\big(\mathcal{E}_k \mathcal{R}_k(\underline{p}, \check{\underline{\zeta}}, 0) - i\mathcal{F}'_k \mathcal{Q}_k(\underline{p}, \check{\underline{\zeta}})\big) \geq \text{Id}\,.$$

5. Summing up, we see that the matrix defined in (5.27) is self adjoint. From (5.30), the condition (5.14) is satisfied provided that $\delta c = 1$. Moreover,

$$\text{Re}\,\big(\mathcal{S}_k(\mathcal{Q}_k + \rho \mathcal{R}_k)\big) = \check{\gamma}\delta \mathcal{D}_k(p,\check{\zeta}) + \rho\delta \mathcal{D}'_k(p,\check{\zeta},\rho),$$

with $\mathcal{D}_k = \text{Diag}(D_k)$ and, at the base point,

$$\mathcal{D}'_k(\underline{p}, \check{\underline{\zeta}}, 0) = \text{Re}\,\big(\mathcal{E}_k \mathcal{R}_k(\underline{p}, \check{\underline{\zeta}}, 0) - i\mathcal{F}'_k \mathcal{Q}_k(\underline{p}, \check{\underline{\zeta}})\big)\,.$$

By (5.32) and (5.34), the matrices $\mathcal{D}_k$ and $\mathcal{D}'_k$ are positive definite on a neighborhood of the base point $(\underline{p}, \check{\underline{\zeta}}, 0)$. This implies that $\mathcal{S}^k$ and $\mathcal{H}_1^k$ satisfy (5.15).

This finishes the construction of symmetrizers for $\mathcal{H}_1^k$ and thus the proof of Theorem 4.1.